\documentclass[12pt,a4paper]{amsart}

\usepackage[left=3cm,right=3cm,top=3.5cm,bottom=3.5cm]{geometry}
\usepackage{amsmath,amsthm,amssymb}
\usepackage[english]{babel}
\usepackage{tikz,pgffor}
\usetikzlibrary {graphs}

\newtheorem{theorem}{Theorem}[section]
\newtheorem{lemma}[theorem]{Lemma}
\newtheorem{prop}[theorem]{Proposition}

\newtheorem{coro}[theorem]{Corollary}

\theoremstyle{definition}
\newtheorem{dfn}[theorem]{Definition}
\newtheorem{rem}[theorem]{Remark}
\newtheorem{exam}[theorem]{Example}

\numberwithin{equation}{section}

\newcommand{\rleft}{\mathopen{}\mathclose\bgroup\left}
\newcommand{\rright}{\aftergroup\egroup\right}

\def\T{{\Bbb T}}
\def\C{{\Bbb C}}

\def\N{{\Bbb N}}
\def\Q{{\Bbb Q}}
\def\R{{\Bbb R}}
\def\P{{\Bbb P}}
\def\Z{{\Bbb Z}}
\def\A{{\Bbb A}}
\def\G{{\Bbb G}}
\def\F{{\Bbb F}}

\newcommand{\ie}{i.\,e.~}

\begin{document}

\title[Classification of smooth toric surfaces]
{On the classification of smooth toric surfaces with exactly one exceptional curve}

\author{Victor  Batyrev}
\address{Fachbereich Mathematik, Universit\"at T\"ubingen, Auf der
Morgenstelle 10, 72076 T\"ubingen, Germany}
\email{batyrev@math.uni-tuebingen.de}
\thanks{}
\dedicatory{Dedicated to Yuri Prokhorov on the occasion of his 60-th birthday}

\begin{abstract}
We classify all smooth projective toric surfaces  $S$ containing exactly one exceptional 
curve. We show that every such surface $S$ is  isomorphic to either $\F_1$ or  a surface $S_r$  defined by a rational number $r \in \Q \setminus \Z$ $(r >1)$. If $a:= [ r]$ then $S_r$ is 
obtained  from the minimal desingularization of the weighted projective plane $\P(1, 2, 2a+1) $ by toric blow-ups whose quantity equals the level  of the rational number  $\{ r \} \in (0,1)$ in the classical Farey tree. Moreover, we show that if $r = b/c$ with coprime $b$ and $c$, then $S_r$ is the minimal desingularization of the weighted projective plane $\P(1, c, b)$. 
We apply   $2$-dimensional regular 
fans $\Sigma_r$ of toric surfaces $S_r$ for constructing  $2$-dimensional colored fans $\Sigma^c$ of minimal horospherical $3$-folds having a  regular $SL(2) \times \G_m$-action. The latter are 
minimal toric $3$-folds $V_r$ classified   by Z. Guan. We establish a direct combinatorial connection between the $3$-dimensional 
fans $\widetilde{\Sigma}^c_r$  of $3$-folds $V_r$ and the $2$-dimensional fans $\Sigma_r$ of surfaces $S_r$.   

\end{abstract}
\maketitle

\thispagestyle{empty}

\section{Introduction}

Smooth projective toric surfaces $S$ are equivariant projective compactifications of  $2$-dimen\-sional 
algebraic torus $\T^2 \cong \G_m^2$ by a circular union of torus invariant smooth projective 
rational curves $C_i \cong \P^1$ $(i \in \Z/n\Z, \, n \geq 3)$  such that every  intersection 
$C_i \cap C_{i+1}$ $(i \in \Z/n\Z)$
is one of  $n$ torus fixed points on $S$.

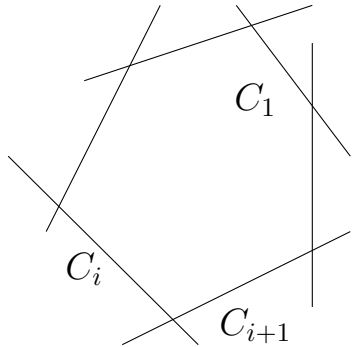
\begin{center}
\begin{tikzpicture}[scale=0.5]
\draw   (4,1) -- (1,5);
\draw  (3,5) -- (-3,3);
\draw  (-1,5) -- (-4,-1);
\draw   (-5,1) -- (0,-4);
\draw  (-2,-4) -- (4,-1);
\draw  (3, -3) -- (3,4);
\node at (-3,-2)  {\large {\bf $C_i$}};
\node at (1.5,-3.5)  {\large {\bf $C_{i+1}$}};
\node at (1.5,2.5)  {\large {\bf $C_1$}};
\end{tikzpicture}
\end{center}

If $M \cong \Z^2$ is the character group of $\T^2$, then every curve $C_i \subset S$ defines  an element $\nu_i \in N :={\rm Hom}(M, \Z)$ which associates to a character $m \in M$, considered as rational function on $S$,  its order along $C_i$. 
The set of rays $\R_{\geq 0} \nu_i$ $(i \in \Z/n\Z)$ in  the plane 
$N_\R := N \otimes \R$ defines a collection of $n$ closed $2$-dimensional 
cones  $\sigma_i:= \R_{\geq 0} \nu_i + \R_{\geq 0}\nu_{i+1}$ ($i \in \Z/n\Z$) with the generators 
$\nu_i, \nu_{i+1}$ forming a $\Z$-basis of $N$. The collection $\Sigma$ consisting of $2$-dimensional cones $\sigma_i$ and  their faces is called  $2$-dimensional {\bf complete regular fan} of toric surface $S$. 
 The fan $\Sigma$,  considered up to unimodular transformation of the lattice $N$ by elements of $GL(2, \Z)$,  
 uniquely determines the toric surface $S$ up to isomorphism. 
 
Since a $GL(2,\Z)$-automorphism of $N$ leaves unchanged linear relations between elements of $N$, it is convenient to combinatorially describe the fan $\Sigma$  
by the associated {\bf weighted circular 
graph}  $\Gamma$ whose $n$ vertices, for simplicity, may be  identified with lattice vectors $\nu_i$ $(i  \in \Z/n\Z)$. The  
{\bf weights $w_i$} attached to vertices $\nu_i$ are determined by $n$ linear relations 
\[ \nu_{i-1} + \nu_{i+1} + w_i \nu_i = 0, \;\; i \in \Z/n\Z.  \]
From the viewpoint of geometry on $S$, the weight $w_i$  equals   
the self-intersection number $\langle C_i, C_i \rangle$ of the rational curve $C_i \subset S$. 
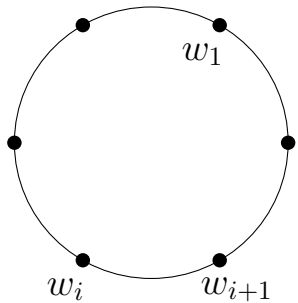
\medskip

\begin{center}
\begin{tikzpicture}[scale = 0.9] 
 \draw (0,0) circle [radius=2cm];
\fill (0:2cm)  circle (3pt);
\fill (60:2cm)  circle (3pt);
\fill (120:2cm)  circle (3pt);
\fill (180:2cm)  circle (3pt);
\fill (240:2cm)  circle (3pt);
\fill (300:2cm)  circle (3pt);
\node at (60:1.5)  {\large {\bf $w_1$}};
\node at (240:2.5)  {\large {\bf $w_i$}};
\node at (300:2.5)  {\large {\bf $w_{i+1}$}};
\end{tikzpicture}
\end{center}

Recall that a smooth projective rational curve $C$ on an arbitrary  smooth 
projective algebraic $S$ is called {\bf exceptional}  if $\langle C, C \rangle = -1$. By theorem of  Castelnuovo, there exists a birational morphism of smooth surfaces $f\, :\, S \to S'$ which contract $C$ to a nonsingular point $p \in S'$ on another smooth surface $S'$. If $S$ is a toric surface, then  every exceptional curve $C$ on $S$ must have empty intersection with 
the affine torus $\T^2 \subset S$. Otherwise such a curve $C$ could be moved by $\T^2$-action and would have non-negative self-intersection number.   Therefore, each 
exceptional curve on toric surface $S$ must be 
one of irreducible $\T^2$-invariant components $C_i$ of $S \setminus \T^2 = \bigcup_{i \in \Z/n\Z} C_i$.  
Moreover, the birational contraction of $C_i \subset S$  is a $\T^2$-equivariant morphism $S \to S'$  which maps $C_i$ to a torus fixed point $p \in S'$ 
on another smooth toric surface $S'$. In this situation, $p$ is a common point of two $\T^2$-invariant 
curves $C_j'$ and $C_{j+1}'$ on $S'$, \ie, $S$ is obtained from $S'$ by the blow-up of the smooth torus fixed point $p = C_j' \cap C_{j+1}'$. 

Our purpose is to classify up to isomorphism all smooth projective toric surfaces $S$ containing exactly one 
exceptional rational curve $C_i$.

In Section 1 we obtain this classification using the  combinatorial  classification of the corresponding weighted circular graphs $\Gamma$ with $n$ vertices such that 
$w_i = \langle C_i, C_i \rangle =-1$ holds true exactly for one $i \in \Z/n\Z$. 
If  $3 \leq n \leq 5$ there exists only one such a weighted circular graph up to isomorphism.  It  corresponds to the  
Hirzebruch surface $\F_1$. If $n=6$ we obtain infinitely many such weighted graphs and denote them $\Gamma(a)$, where  $a$ is a positive integer.  Furthermore,  we show that  if  $n \geq 6$, then for each  weighted circular graph $\Gamma$ with $n$ vertices and with exactly one weight $-1$ there exists a positive integer $a$ such that $\Gamma$ is obtained from the  graph $\Gamma(a)$ by a unique  finite sequence of  $n-6$ combinatorial blow-ups.  For every 
 integer $a \geq1$, we consider the  rooted oriented tree $\mathcal{OT}(a)$ whose vertices are weighted circular graphs $\Gamma$ with $n \geq 6$ vertices containing exactly one weight $-1$, whose oriented edges $(\Gamma, \Gamma')$ correspond to combinatorial blow-ups $\Gamma \to \Gamma'$, and whose 
root node is $\Gamma(a)$.  We show that  $\mathcal{OT}(a)$ is a complete  infinite rooted 
binary tree  with the  root node $\Gamma(a)$.  
 It is convenient to identify $\mathcal{OT}(a)$ with  the classical Farey tree  $F$ whose nodes bijectively parametrize all rational numbers $\delta \in (0,1)$. In this identification, the  root node $\frac{1}{2} \in (0,1)$ of the Farey tree $F$ corresponds to the weighted circular graph $\Gamma(a)$. This allows to enumerate  nodes  $\Gamma$ in the oriented  tree $\mathcal{OT}(a)$ by rational numbers $\delta \in (0,1)$. We denote these nodes  by $\Gamma_\delta(a)$ and provided the  alternative notation $\Gamma_{\frac{1}{2}}(a)$ for the root node
$\Gamma(a)$. In order to shortly cover  in a unified way all nodes $\Gamma$ in  all trees $\mathcal{OT}(a)$ $(a \geq 1)$  we use the rational number $r:= a + \delta \in \Q \setminus \Z$ $(r > 1)$ and denote 
 the weighted circular graphs $\Gamma_\delta(a) \in \mathcal{OT}(a)$ simply by $\Gamma_r$. 

In Section 2 we  study the  geometry of the smooth toric surfaces $S_r  $ corresponding to the 
weighted circular graphs $\Gamma_r$ obtained from the above classification. First, we use the minimal combinatorial desingularization of $2$-dimensional toric singularities described by 
$2$-dimensional rational cones $\sigma \subset N_\R$. 
Next, we show that each  graph $\Gamma(a) = \Gamma_{a +\frac{1}{2}}$ describes a smooth toric surface $S_{a+\frac{1}{2}}$ which is the minimal desingularization $\widehat{\P}(1,2,2a+1)$ of the weighted projective plane $\P(1, 2, 2a+1)$. 
Finally, we show that if   $r = b/c \in \Q \setminus \Z (r>1)$ with coprime $b, c \in \N$, then    
$S_r$ is isomorphic to  the minimal desingularization of the weighed projective plane $\P(1, c, b)$.  

Section 3 is devoted to an application of our classification. We use $2$-dimensional complete regular fans $\Sigma_r$ defining toric surfaces 
$S_r$ for constructing $2$-dimensional colored fans $\Sigma_r^c$ which define minimal horospherical $3$-folds $V_r$ having  
a regular $SL(2) \times \G_m$-action. These $3$-folds are  horospherical compactifications of the $3$-dimensional homogeneous space 
$V:=(SL(2)/U ) \times \G_m$,  where $U \subset SL(2)$ is a maximal unipotent subgroup. The horospherical $3$-folds 
$V_r$ are defined by the colored $2$-dimesional fans $\Sigma_r^c$ obtained from $\Sigma_r$ by attaching 
to the unique ray $\R_{\geq 0} \nu_i$ in $\Sigma_r$ with the  weight $w_i =-1$ the unique color of $V$. Using the spherical Mori theory \cite{Br93,Br94}, we show that  the horospherical varieties $V_r$ associated with 
$\Sigma_r^c$ are minimal, \ie, they do not admit any non-trivial  equivariant Mori contraction to another  smooth projective $3$-dimensional horospherical $3$-fold. On the other hand,  $V_r$ are minimal smooth toric $3$-folds classified by Z. Guan \cite{Gu00}. We show how  the $3$-dimensional fan of  the toric $3$-fold  $V_r$ can be easily obtained from the $2$-dimensional toric  fan of $S_r$. 

In Appendix we shorly summarize the Luna-Vust theory of spherical varieties and Mori theory of horospherical varieties of rank $2$ and explain how the 
classification of smooth toric surfaces with exactly one exceptional curve can be applied   to constructing  minimal smooth projective $G$-equivariant compactifications of  horospherical homogeous spaces $G/H$ of rank $2$ having exactly one color described by 
a primitive lattice vector $\nu$ in the valuation lattice $N \cong \Z^2$. 

{\bf Acknowledgements.} The author is grateful for discussions of the presented topic with V. Alexeev, Yu. Prokhorov and F. Haddad. 

\section{The classification of weighted circular graphs}

 Throughout this paper we use the bijective correspondence between the following two sets \cite{Oda85}: 
\begin{itemize}
\item 
Smooth projective toric surfaces $S$ considered up to isomorphism.
\item  
Complete regular $2$-dimensional fans $\Sigma$ in $N_\R$ considered 
up to unimodular isomorphism. 
\end{itemize}

In our classification we also need the following statements explained above:

\begin{prop} \label{isomorphism} \cite[Cor.~1.29]{Oda85} Two complete regular $2$-dimenional fans $\Sigma$ and $\Sigma'$ are unimodular equivalent to each other  if and only if the corresponding weighted circular graphs $\Gamma$ and $\Gamma'$ are isomorphic as weighted graphs. 
\end{prop} 

\begin{prop} \label{-1-curve} A smooth projective toric surface $S$ contains exactly one exceptional $(-1)$-curve if and only if 
the corresponding weighted circular graph $\Gamma$ contains exactly one weight $-1$.
\end{prop}

\begin{dfn}
We say that a weighted circular graph $\Gamma'$ with $n+1$ vertices 
is obtained from the weighted circular graph $\Gamma$ with $n$ vertices
by a {\bf combinatorial blow-up} ( or simply by a {\bf blow-up}) if there exists an edge $[\nu_i, \nu_{i+1}]$  of $\Gamma$ such that  $\Gamma'$ corresponds to the  fan $\Sigma'$
 obtained from a complete regular fan $\Sigma$ corresponding to $\Gamma$ by the subdivision of the cone 
$\R_{\geq 0}\nu_i + \R_{\geq 0} \nu_{i+1}$ into union of two regular cones 
\[ \R_{\geq 0}\nu_i + \R_{\geq 0} \nu' \;\; {\rm and } \;\;   \R_{\geq 0}\nu' + \R_{\geq 0} \nu_{i+1}, \]
where $\nu':= \nu_i + \nu_{i+1}$. 
\end{dfn}

\begin{rem}\label{change-weights}
If a weighted circular graph $\Gamma'$ is obtained from $\Gamma$ by a combinatorial blow-up as above, then all weights
of $\Gamma$, except $w_i$ and $w_{i+1}$, remains unchanged and one has 
$w_i'= w_i -1$, $w'_{i+1} = w_{i+1}-1$, 
$w':= -1$, \ie,  we obtain the following transformation of weighted circular graphs $\Gamma \to \Gamma'$:

\begin{center}
\begin{tikzpicture}[scale = 0.8]
 \draw (0,0) circle [radius=2cm];
\fill (0:2cm)  circle (3pt);
\fill (60:2cm)  circle (3pt);
\fill (120:2cm)  circle (3pt);
\fill (180:2cm)  circle (3pt);
\fill (225:2cm)  circle (3pt);
\fill (315:2cm)  circle (3pt);
\fill (270:2cm)  circle (3pt);
\node at (60:1.5)  {\large {\bf $w_1$}};
\node at (225:2.7)  { {\bf $w_i -1$}};
\node at (315:2.7)  { {\bf $w_{i+1} -1$}};
\node at (270:2.5)  { {\bf $-1$}};

\begin{scope}[xshift= -7cm] 
 \draw (0,0) circle [radius=2cm];
\fill (0:2cm)  circle (3pt);
\fill (60:2cm)  circle (3pt);
\fill (120:2cm)  circle (3pt);
\fill (180:2cm)  circle (3pt);
\fill (240:2cm)  circle (3pt);
\fill (300:2cm)  circle (3pt);
\node at (60:1.5)  {\large {\bf $w_1$}};
\node at (240:2.5)  {\large {\bf $w_i$}};
\node at (300:2.5)  {\large {\bf $w_{i+1}$}};
\draw [ - > ]  (3,0) -- (4,0);
\end{scope}
\end{tikzpicture}
\end{center}
\end{rem}

\begin{theorem} \cite{Oda85}
If $n=3$, then the weighted circular graph is isomorphic to 
 
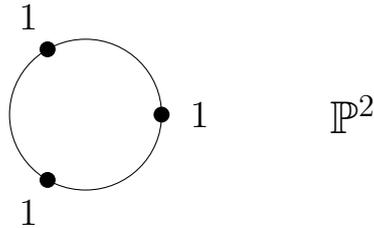
\begin{center}
\begin{tikzpicture}
 \draw (0,0) circle [radius=1cm];
\fill (0:1cm)  circle (3pt);
\fill (120:1cm)  circle (3pt);
\fill (240:1cm)  circle (3pt);

\node at (0:1.5)  {\large {\bf $1$}};
\node at (120:1.5)  {\large {\bf $1$}};
\node at (240:1.5)  {\large {\bf $1$}};

\begin{scope}[xshift = 100]
\node at (0,0) {\Large $\P^2$};
\end{scope}

\end{tikzpicture}
\end{center}

If $n=4$, the the weighted circular graph is isomorphic to 
 
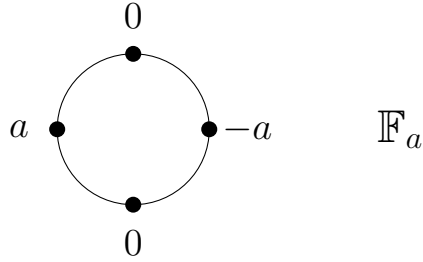
\begin{center}
\begin{tikzpicture}
 \draw (0,0) circle [radius=1cm];
\fill (0:1cm)  circle (3pt);
\fill (90:1cm)  circle (3pt);
\fill (180:1cm)  circle (3pt);
\fill (270:1cm)  circle (3pt);

\node at (0:1.5)  {\large {\bf $-a$}};
\node at (90:1.5)  {\large {\bf $0$}};
\node at (180:1.5)  {\large {\bf $a$}};
\node at (270:1.5)  {\large {\bf $0$}};

\begin{scope}[xshift = 100]
\node at (0,0) {\Large $\F_a$};
\end{scope}
\end{tikzpicture}
\end{center}

Every weighted circular graph with $n \geq 5$ vertices is obtained from a weighted graph with 
$4$ vertices by a sequence of $n-4$ combinatorial  blow-ups. 
\end{theorem} 
 
\begin{coro}\label{n=5}
Every weighted circular graph $\Gamma$ with $n=5$ vertices is isomorphic to the graph 
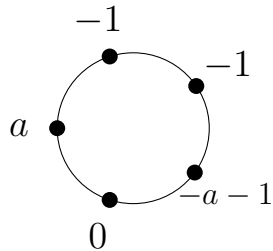
\begin{center}
\begin{tikzpicture}
 \draw (0,0) circle [radius=1cm];
\fill (34:1cm)  circle (3pt);
\fill (108:1cm)  circle (3pt);
\fill (180:1cm)  circle (3pt);
\fill (252:1cm)  circle (3pt);
\fill (324:1cm)  circle (3pt);

\node at (34:1.5)  {\large {\bf $-1$}};
\node at (108:1.5)  {\large {\bf $-1$}};
\node at (180:1.5)  {\large {\bf $a$}};
\node at (252:1.5)  {\large {\bf $0$}};
\node at (324:1.5)  { {\bf $-a-1$}};
\end{tikzpicture}
\end{center}
for some integer $a$. In particular, $\Gamma$ always contains at least two vertices with weights $-1$. 
\end{coro}
 
Applying a combinatorial blow-up to the graph in \ref{n=5}, we obtain the following $3$ possible forms of weighted circular graph $\Gamma$ with $6$ vertices:

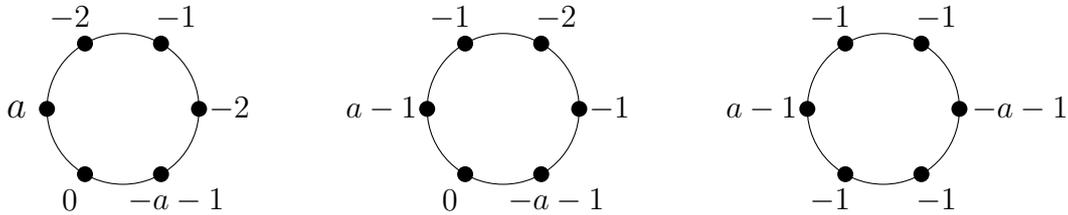
\begin{center}

\begin{tikzpicture}
 \draw (0,0) circle [radius=1cm];
\fill (60:1cm)  circle (3pt);
\fill (0:1cm)  circle (3pt);
\fill (120:1cm)  circle (3pt);
\fill (180:1cm)  circle (3pt);
\fill (240:1cm)  circle (3pt);
\fill (300:1cm)  circle (3pt);

\node at (60:1.4)  {{\bf $-1$}};
\node at (120:1.4)  { {\bf $-2$}};
\node at (180:1.4)  {\large {\bf $a$}};
\node at (240:1.4)  { {\bf $0$}};
\node at (300:1.4)  { {\bf $-a-1$}};
\node at (0:1.4)  {{\bf $-2$}};

\begin{scope}[xshift = 5 cm]
 \draw (0,0) circle [radius=1cm];
\fill (60:1cm)  circle (3pt);
\fill (0:1cm)  circle (3pt);
\fill (120:1cm)  circle (3pt);
\fill (180:1cm)  circle (3pt);
\fill (240:1cm)  circle (3pt);
\fill (300:1cm)  circle (3pt);

\node at (60:1.4)  {{\bf $-2$}};
\node at (120:1.4)  { {\bf $-1$}};
\node at (180:1.6)  { {\bf $a-1$}};
\node at (240:1.4)  { {\bf $0$}};
\node at (300:1.4)  { {\bf $-a-1$}};
\node at (0:1.4)  {{\bf $-1$}};
\end{scope}

\begin{scope}[xshift = +10 cm]
 \draw (0,0) circle [radius=1cm];
\fill (60:1cm)  circle (3pt);
\fill (0:1cm)  circle (3pt);
\fill (120:1cm)  circle (3pt);
\fill (180:1cm)  circle (3pt);
\fill (240:1cm)  circle (3pt);
\fill (300:1cm)  circle (3pt);

\node at (60:1.4)  {{\bf $-1$}};
\node at (120:1.4)  { {\bf $-1$}};
\node at (180:1.6)  {{\bf $a-1$}};
\node at (240:1.4)  { {\bf $-1$}};
\node at (300:1.4)  { {\bf $-1$}};
\node at (0:1.8)  {{\bf $-a-1$}};
\end{scope}
\end{tikzpicture}

\end{center}

 \begin{coro}\label{n=6}
Every weighted circular graph $\Gamma$ with $n=6$ vertices and containing exactly one weight -1 is isomorphic to the following graph  $\Gamma(a)$ $(a \geq 1)$:
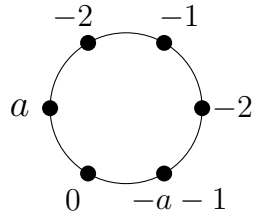
\begin{center}
\begin{tikzpicture}
 \draw (0,0) circle [radius=1cm];
\fill (60:1cm)  circle (3pt);
\fill (0:1cm)  circle (3pt);
\fill (120:1cm)  circle (3pt);
\fill (180:1cm)  circle (3pt);
\fill (240:1cm)  circle (3pt);
\fill (300:1cm)  circle (3pt);

\node at (60:1.4)  {{\bf $-1$}};
\node at (120:1.4)  { {\bf $-2$}};
\node at (180:1.4)  {\large {\bf $a$}};
\node at (240:1.4)  { {\bf $0$}};
\node at (300:1.4)  { {\bf $-a-1$}};
\node at (0:1.4)  {{\bf $-2$}};
\end{tikzpicture}
\end{center} 
 
\end{coro}

\noindent 
{\bf Proof.} The second and third graphs contain already at least two weights $(-1)$. The first graph contains
exactly one weight $(-1)$ if $a \not\in \{ 0, -1 \}$. If we denote the first graph by $\Gamma(a)$, then two graphs 
$\Gamma(a)$ and $\Gamma(-a-1)$ are isomorphic. Therefore, it is enough to consider only the graphs $\Gamma(a)$ for $a \geq 1$. \hfill $\Box$
 \medskip
 
In our combinatorial classification we use  the following simple statement:

\begin{lemma} \label{two-weights}
Let $\Gamma$ be a weighted circular graph $\Gamma$ containing  at least 
two weights $-1$ belonging to non-neighbour vertices. Then 
a weighted circular graph $\Gamma'$ obtained from 
$\Gamma$ by a single combinatorial blow-up also contains two  
weights $-1$ belonging to non-neighbour vertices.
\label{two-1} 
\end{lemma}

\noindent
{\bf Proof.} Let $\nu':= \nu_i + \nu_{i+1}$. It follows from our assumption
about $\Gamma$ that there exists a vertex $\nu_j$  $(j \neq i, j \neq i+1)$
such that $w_j =-1$, because  the vertices $v_i$ and $v_{i+1}$ 
are neighbour to each other. By \ref{change-weights},  the weight $w_j=-1$ remains   
unchanged after the blow-up 
$\Gamma \to \Gamma'$ and we get one more weight $w'=-1$ corresponding to the new vertex $\nu'$. 
Obviously, $\nu'$ and $\nu_j$ are not neighbour to each other. 
\hfill $\Box$

\begin{theorem} \label{chain}
Let $\Gamma$ be  a weighted circular graph with $n \geq 6$ vertices containing  
exactly  one vertex with the weight $-1$.  Then there exists a positive number $a \in \Z$ such that 
either $n =6$ and 
$\Gamma$ is isomorphic  
to $\Gamma(a)$, or $\Gamma$ is obtained from $\Gamma(a)$ by 
a uniquely determined sequence of combinatorial blow-ups 
$$\Gamma(a) = \Gamma^0 \to \Gamma^1 \to \cdots \to \Gamma^{n-7} \to 
\Gamma^{n-6} = \Gamma$$
such that every weighted circular graph $\Gamma^l$ $( 0 \leq l \leq n-6)$ 
contains  exactly $3$ vertices with the weights at least $ -1$: 
one vertex $\nu_i$ with the weight 
$-1$, and two neighbour to each other vertices $\nu_j, \nu_{j+1}$  with weights 
$\{ w_j, w_{j+1} \} = \{ 0, a \}$.
Moreover, the unique vertex of $\Gamma^{l}$ with weight $-1$ is always not a neighbour of $\nu_j$ or $\nu_{j+1}$. 
\label{sequence}
\end{theorem} 

\noindent
{\bf Proof.} We prove the statement by induction on $n \geq 6$. 
The case $n=6$ obviously follows  from Corollary \ref{n=6}. 

Consider the next case $n =7$.  Let 
$\Gamma_{-1}$ be the  weighted circular graph with $6$ vertices 
obtained from $\Gamma$ 
by contracting the single vertex $\nu_1$ with the weight $-1$. 
By Lemma \ref{two-1}, $\Gamma'$ can not contain $2$ non-neighbour 
vertices with weights $-1$. By \ref{n=6},  the weighted circular graph 
$\Gamma_{-1}$ must be isomorphic to 
$\Gamma(a)$ for some $a \geq 1$, and all statement of Theorem \ref{sequence} 
hold true.     

Now we assume that $n \geq 8$. Take again 
$\Gamma_{-1}$ to be the  weighted circular graph with $n-1$ vertices 
obtained from $\Gamma$ 
by contracting the single vertex $\nu_i$  with the weight $w_i=-1$. To be able to apply induction hypothesis, we must show  
that $\Gamma_{-1}$ also contains exactly one vertex with weight $-1$. 
By  \ref{change-weights} and \ref{two-1}, $\Gamma_{-1}$ can not contain $2$ non-neighbour 
vertices with the weight $-1$. Therefore, either $\Gamma_{-1}$ contains
exactly one vertex from $\{\nu_{i+1}, \nu_{i-1}\}$ 
with the weight $-1$, or  $\Gamma_{-1}$ contains
both neighbour vertices $\nu_{i+1}, \nu_{i-1}$ with the weight $-1$. 
 In the first case, we are done.  
 
 Let us consider the second case. We  again 
blow-down a $(-1)$-vertex of $\Gamma_{-1}$ 
(say $\nu_{i+1}$) and obtain the graph $\Gamma_{-2}$ 
with $n-2$ vertices. 

By \ref{change-weights}, $\Gamma_{-2}$ can have exactly one vertex $\nu_{i+2}$ with the weight $-1$, 
because the weight of the vertex $\nu_{i-1}$ in $\Gamma_{-2}$ becomes $0$. 
Since $n \geq 8$, we can apply the induction assumption to $\Gamma_{-2}$. This 
yields a contradiction because, by induction hypothesis,  the single possible vertex  with the 
weight $-1$ in $\Gamma_{-2}$ (i.e.,  $\nu_{i_+2}$) can not be  neighbour of  
$\nu_{i-1} \in \Gamma_{-2}$ with the weight $0$. 
So the second case is impossible and $\Gamma_{-1}$ must always contain exactly one  weight $-1$.

Now we can appy the induction hypothesis to $\Gamma_{-1}$ and obtain a sequence of combinatorial blow-ups
$$\Gamma(a) = \Gamma^0 \to \Gamma^1 \to \cdots \to \Gamma^{n-8} \to 
\Gamma^{n-7} = \Gamma_{-1}$$
which satisfies all properties of Theorem \ref{sequence}. Since every graph $\Gamma^k$  in this sequence has exactly 
one vertex with weight $-1$, we obtain a unique sequence of combinatorial blow-ups connecting 
$\Gamma$ with $\Gamma(a)$ for some positive integer $a$. 
\hfill $\Box$ 

\begin{coro}
Let $\mathcal{OT}$ be the oriented tree whose nodes is the  set of all weighted circular graphs $\Gamma$ with $n \geq 6$ vertices having exactly one weight $-1$ and whose oriented edges are pairs $(\Gamma, \Gamma')$, where 
$\Gamma'$ is a combinatorial blow-up of $\Gamma$. Then $\mathcal{OT}$ is a disjoint union of infinite complete binary trees 
$\mathcal{OT}(a)$ parametrized by their root nodes $\Gamma(a)$, $a \in \Z_{>0}$.  
\end{coro}

\noindent
{\bf Proof.} Let $\Gamma$ be a  weighted circular graphs $\Gamma$ with $n \geq 6$ vertices having exactly one 
vertex $\nu_i$ with the weight $w_i=-1$.  By \ref{change-weights},   there exist exactly two possibilities for combinatorial blow-ups $\Gamma \to \Gamma'$ such that $\Gamma'$ also has exactly one weight $-1$: either one blows up the edge $[\nu_{i+1}, \nu_i]$, or the edge 
$[\nu_i, \nu_{i-1}]$.  So starting from the root node $\Gamma(a) \in \mathcal{OT}(a)$ we have exactly $2^k$ possibilities for creating  a path of length $k$ in $\mathcal{OT}(a)$. 
\hfill $\Box$

\begin{exam} We  illustrate the process of subsequent blow-ups using the regular $2$-dimensional cone $\sigma \subset \R^2$ spanned by  
$e_1:=(1,0)$ and $e_2:=(1,1)$. The combinatorial blow-up subdivides $\sigma$ into the union 
of two smaller regular $2$-dimensional cones: 
$$\sigma = \sigma_1 \cup \sigma_2 = \left( \R_{\geq 0} (1,0) + \R_{\geq 0} (1,2)\right) \bigcup \left( \R_{\geq 0} (1,2) + \R_{\geq 0} (1,1)\right) $$ 
having the common 
middle ray $\R_{\geq 0} (e_1 + e_2) = \R_{\geq 0}(1,2)$. Next we blow up each of two cones $\sigma_1$, $\sigma_2$ 
to obtain a subdivision of $\sigma$ in four smaller cones and repeat this process with each obtained subcone again. At each $k$-th step we subdivide $\sigma$ into $2^k$ smaller regular subconed using new rays with the rational slope $\delta$ in the interval $(0,1)$. 

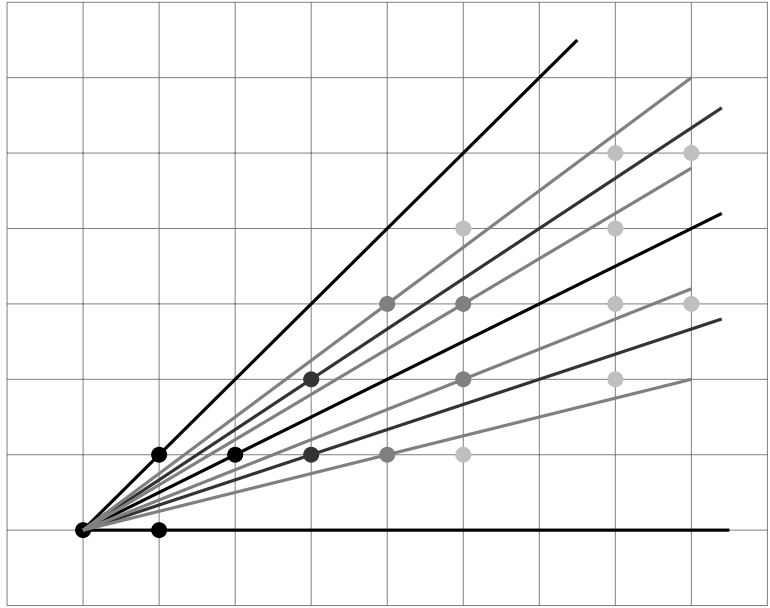
\begin{center}

\begin{tikzpicture}
\draw[step=1cm,gray,very thin]
(-1,-1) grid (9,7) ;

\fill (0,0) circle (3pt);
\fill (1,1) circle (3pt);
\fill (1,0) circle (3pt);
\fill  (2,1) circle (3pt);
\fill [color =black!80!white] (3,2) circle (3pt);
\fill [color =black!80!white] (3,1) circle (3pt);

\fill [color =black!50!white] (4,1) circle (3pt);
\fill [color =black!50!white] (5,2) circle (3pt);
\fill [color =black!50!white] (5,3) circle (3pt);
\fill [color =black!50!white] (4,3) circle (3pt);

\fill [color =black!25!white] (5,1) circle (3pt);
\fill [color =black!25!white] (7,2) circle (3pt);
\fill [color =black!25!white] (8,3) circle (3pt);
\fill [color =black!25!white] (7,3) circle (3pt);
\fill [color =black!25!white] (7,4) circle (3pt);
\fill [color =black!25!white] (8,5) circle (3pt);
\fill [color =black!25!white] (7,5) circle (3pt);
\fill [color =black!25!white] (5,4) circle (3pt);

\draw [solid, very thick] (0,0) -- (6.5,6.5);
\draw  [solid, very thick]  (0,0) -- ( 8.5, 0);
\draw [ solid, very thick ] (0,0) -- (8.4, 4.2);
\draw [color =black!80!white, solid, very thick] (0,0) -- (8.4, 2.8);
\draw [color =black!80!white, solid, very thick] (0,0) -- (8.4, 5.6); 

\draw [color =black!50!white , solid, very thick] (0,0) -- (8,2);
\draw [color =black!50!white , solid, very thick] (0,0) -- (8,3.2);
\draw [color =black!50!white , solid, very thick] (0,0) -- (8,4.8);
\draw [color =black!50!white , solid, very thick] (0,0) -- (8,6);
\end{tikzpicture}

\end{center}

The combinatorial blow-up of a regular 
subcone $\widetilde{\sigma} \subset \sigma$ generated by two primitive lattice vectors $u_1:=(a,b)$, $u_2:= (c,d)$ satisfying 
the condition  
$ad-bc =1$ adds the new lattice vector $u:= (a+c, b+d)$ whose slope $\frac{a+c}{b+d}$ is called {\bf Farey mediant sum} of 
the slopes $\frac{a}{b}$ and $\frac{c}{d}$:
\[ \frac{a}{b} \oplus \frac{c}{d} = \frac{a+c}{b+d}. \]
The classical {\bf Farey  tree} is a binary tree with the root $\frac{1}{2} \in (0,1)$ containing all rational numbers from $(0,1)$ in some ordered way, constucted hierarchically, level by level, using the Farey mediant sum. 

 
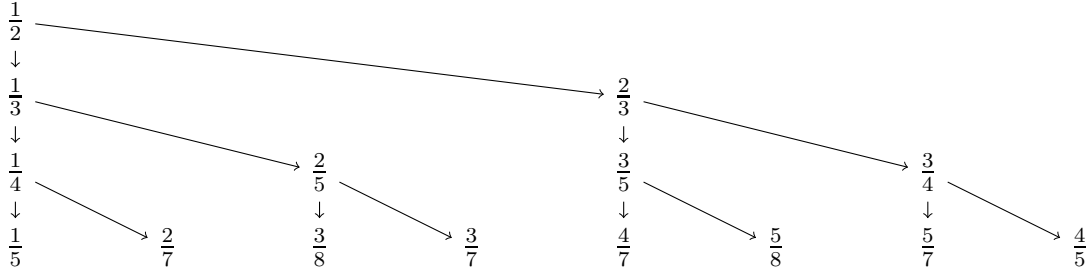
\tikz
 \graph [grow down,
          branch right=2cm] {
"$\frac{1}{2}$" -> {"$\frac{1}{3}$" ->{"$\frac{1}{4}$" 
->{"$\frac{1}{5}$","$\frac{2}{7}$"} ,"$\frac{2}{5}$"  ->{"$\frac{3}{8}$","$\frac{3}{7}$"}  }, "$\frac{2}{3}$" 
->{"$\frac{3}{5}$"  ->{"$\frac{4}{7}$","$\frac{5}{8}$"}  ,"$\frac{3}{4}$" 
 ->{"$\frac{5}{7}$","$\frac{4}{5}$"}  }}
};
\end{exam}

\begin{dfn}
Now we associate with every weighted circular graph $\Gamma \in \mathcal{OT}$ a rational number $r \in \Q \setminus \Z$ $(r >1)$ defined as $r := a + \delta$, where the integer $a$ is defined by the condition $\Gamma \in \mathcal{OT}(a)$ and the rational number $\delta \in (0,1)$ describes  the position of the node $\Gamma$ in the classical Farey after its natural identification 
with the binary tree   $\mathcal{OT}(a)$ with the root node $\Gamma(a)$. Now  we denote the weighted circular graph $\Gamma \in \mathcal{OT}$ simply by $\Gamma_r$. The corresponding $2$-dimensional  complete regular fan will be denoted by $\Sigma_r$, and the corresponding smooth projective toric surface by $S_r$. In these new notations, the root node  $\Gamma(a) \in 
\mathcal{OT}(a)$ is $\Gamma_{a + \frac{1}{2}}$. 
\end{dfn} 

\section{Minimal desingularizations of weighted projective planes}

Now we want to give a more geometric description of all smooth toric surfaces $S_r$. For this purpose, up to unimodular transformation, we choose 
two neighbour vertices $\nu_{j}$, $\nu_{j+1}$ from Theorem \ref{chain} as the following lattice vectors in $N = \Z^2$:
\[  \nu_j:= (-1,0), \nu_{j+1} := (0, -1), \;\; w_j := a \geq 1, \,\, w_{j+1} := 0. \]

\begin{exam} \label{1-2-3} Consider examples of $\Sigma_r$ for $r = a + \delta$ with 
$\delta = \frac{1}{2}$ and $a \in \{ 1, 2, 3 \}$:

\begin{center}
\begin{tikzpicture}
\draw[step=1cm,gray,very thin]
(-2,-2) grid (4,3) ;

\node at (0, -1.5) {$0$};
\node at (-1.5, 0) {$ 1$};

\node at (1,1.5) {$-2$};
\node at (3, 2.5) {$-1$};
\node at (1, -0.5) {$-2$};
\node at (2.5, 1) {$-2$};

\fill (0,0) circle (3pt);
\fill (-1,0) circle (3pt);
\fill (0,-1) circle (3pt);
\fill (3,2) circle (3pt);
\fill (1,1) circle (3pt);
\fill (2,1) circle (3pt);
\fill (1,0) circle (3pt);

\draw  [solid, very thick]  (0,0) -- ( -1, 0);
\draw  [solid, very thick]  (0,0) -- ( 0, -1);
\draw  [solid, very thick] (0,0) -- (1,0); 

\draw [solid, very thick] (0,0) -- (1,1);
\draw [solid, very thick] (0,0) -- (2,1);
\draw [solid, very thick] (0,0) -- (3,2);

\node at (2.5,-0.5) {$a=1$};

\draw [solid, very thin]  (-1,0) -- (0,-1) -- (1,0)  -- (3,2) -- (-1,0) ;

\begin{scope}[xshift = 5cm]
\draw[step=1cm,gray,very thin]
(-2,-2) grid (6,3) ;

\fill (0,0) circle (3pt);
\fill (-1,0) circle (3pt);
\fill (0,-1) circle (3pt);

\fill (5,2) circle (3pt);
\fill (2,1) circle (3pt);
\fill (3,1) circle (3pt);

\fill (1,0) circle (3pt);

\draw  [solid, very thick]  (0,0) -- ( -1, 0);
\draw  [solid, very thick]  (0,0) -- ( 0, -1);
\draw  [solid, very thick] (0,0) -- (1,0); 

\draw [solid, very thick] (0,0) -- (2,1);
\draw [solid, very thick] (0,0) -- (3,1);
\draw [solid, very thick] (0,0) -- (5,2);

\node at (0, -1.5) {$0$};
\node at (-1.5, 0) {$2$};

\node at (2,1.5) {$-2$};
\node at (5, 2.5) {$-1$};
\node at (1, -0.5) {$-3$};
\node at (3.5, 1) {$-2$};

\node at (3,-0.5) {$a=3$};

\draw [solid, very thin]  (-1,0) -- (0,-1) -- (1,0)  -- (5,2) -- (-1,0) ;
\end{scope}
\end{tikzpicture}

\end{center}

\begin{center}
\begin{tikzpicture}
\draw[step=1cm,gray,very thin]
(-2,-2) grid (8,3) ;

\fill (0,0) circle (3pt);
\fill (-1,0) circle (3pt);
\fill (0,-1) circle (3pt);

\fill (7,2) circle (3pt);
\fill (3,1) circle (3pt);
\fill (4,1) circle (3pt);

\fill (1,0) circle (3pt);

\draw  [solid, very thick]  (0,0) -- ( -1, 0);
\draw  [solid, very thick]  (0,0) -- ( 0, -1);
\draw  [solid, very thick] (0,0) -- (1,0); 

\draw [solid, very thick] (0,0) -- (3,1);
\draw [solid, very thick] (0,0) -- (4,1);
\draw [solid, very thick] (0,0) -- (7,2);

\node at (0, -1.5) {$0$};
\node at (-1.5, 0) {$ 3$};

\node at (3,1.5) {$-2$};
\node at (7, 2.5) {$-1$};
\node at (1, -0.5) {$-4$};
\node at (4.5, 1) {$-2$};

\node at (5,-0.5) {$a=3$};

\draw [solid, very thin]  (-1,0) -- (0,-1) -- (1,0)  -- (7,2) -- (-1,0) ;
\end{tikzpicture}
\end{center}

\end{exam}

\begin{rem} \cite{Oda85} \label{minimal}
Recall that every affine toric surface singularity is a $2$-dimensional cyclic quotient singularity $X \cong \A^2/\mu_d$ of order $d$ which is described by a $2$-dimensional 
rational cone $\sigma = \rho_1\R_{\geq 0} + \rho_2 \R_{\geq 0} \subset \N_\R \cong \R^2$ where $\rho_1 =(a,b)$ and $\rho_2 = (b,c)$ are primitive lattice vectors $N$ satisfying the condition  $|ad-bc| =d$. A minimal toric desingularization $f\,:\, \widehat{X} \to X$ is determined by a refinement of $\sigma$ using $k$ primitive lattice vectors $v_1, \ldots , v_k $ in the interior of $\sigma$ such that $k+2$ lattice points $\rho_1, v_1, \ldots, v_k, \rho_2$ are exactly all lattice points on the compact boundary of the convex hull ${\rm Conv}(\sigma \cap (N \setminus \{(0,0)\} \}$.  If we put $v_0:= \rho_1$ and 
$v_{k+1}:= \rho_2$ then every pair $\{ v_{i-1}, v_{i} \}$ $( 1 \leq i \leq k+1)$ form a $\Z$-basis of $N$ and one obtains  $k$ integers 
$a_1, \ldots, a_k$ $(a_i \leq -2)$ satisfying the conditions 
\[  v_{i-1} + v_{i+1} + a_i v_i =0,  \;\;;\; \forall \; i  \in \{1, \ldots, k\}, \]
$$ \det  \begin{pmatrix} a_1 & 1 & 0 & \cdots & 0 \\
1 & a_2 & 1 & \cdots, & 0  \\
0 & 1 & a_3 & \cdots & 0  \\
\cdot & \cdot & \cdot & \cdots & \cdot \\
0 & 0 &  0 & \cdots & a_k 
\end{pmatrix} = (-1)^k d.  $$
The exceptional fiber $f^{-1}(0)$ in the minimal desingularization is a chain of $k$ smooth projective rational cuves 
$C_1, \ldots, C_k$ and $a_i = \langle C_i, C_i \rangle$, $i \in \{1, \ldots, k\}$.
\end{rem}

Examples in \ref{1-2-3} illustrate the following description  of the toric surfaces $S_{a+ \frac{1}{2}}$:

\begin{prop}\label{1/2}
Every toric surface $S_{a+ \frac{1}{2}}$ is isomorphic to the minimal desingularization $\widehat{\P}(1,2,2a+1)$ of the singular weighted projective plane $\P(1,2,2a+1)$.  
\end{prop}

\noindent 
{\bf Proof.} It follows from the relation $\nu_{j-1} + a \nu_j + \nu_{j+1} =0$ that $\nu_{j-1} = (1, a)$. 
Then we have the relation $\nu_{j-2} + (-2)\nu_{j-1} + \nu_j = 0$ which implies that $\nu_{j-2} = (2, 2a+1)$. 
Note that the weight $w_{j-2}$ equals $-1$. The lattice vectors $\nu_{j-2}$, $\nu_j$ and $\nu_{j+1}$ satisfy 
the linear relation 
\[ \nu_{j-2} + (2a+1) \nu_j + 2\nu_{j+1} = 0. \]
So they span the fan of the weighted projective plane $\P(1,2, 2a+1)$. This weighted projective plane contains two 
singular points which are cyclic quotient singularities of order $2$ and $2a+1$. The minimal resolution of the first singularity consists of one smooth rational curve $C$ with the self-intersection number $-2$. The minimal resolution of the second singularity consists of the union of two intersecting smooth rational curves $C'$ and $C''$ with the self-intersection numbers $-2$ and $-a-1$. 
\hfill $\Box$  

The last statement can be generalized for arbitrary toric surface $S_r$.

\begin{theorem}
Let $r \in \Q \setminus \Z$ $(r>1)$. Choose the lattice generators 
$\nu_1$, $\nu_j$ and $\nu_{j+1}$  in the fan $\Gamma_r$ such that $w_1 =-1$, $w_j = a >0$ and $w_{j+1} = 0$. Then 
the integral linear relation among $\nu_1$, $\nu_j$ and $\nu_{j+1}$ has 
the form 
\[ \nu_1 + b \nu_j + c \nu_{j+1} = 0, \;\; b>c >1, \gcd (b,c) = 1,  \] 
where $r = b/c$. 
Moreover, the toric surface $S_r$ is isomorphic to the minimal desingularization  $\widehat{\P}(1, c, b)$ of the  singular weighted projective 
plane  ${\P}(1, c, b)$. 
\end{theorem}

\noindent
{\bf Proof.} Since $\nu_j$ and $\nu_{j+1}$ form a $\Z$-basis of $N$ we can write $\nu_1 +  b \nu_j + c \nu_{j+1} = 0$ with 
integer coefficients $b$ and $c$. Since all weights $w_2, \ldots, w_{j-1}$ are not greater than $-2$, by \ref{minimal}, the lattice vectors $\nu_2, \ldots, \nu_{j-1}$ define the minimal regular refinement of the cone  
$\R_{\geq 0} \nu_1 + \R_{\geq 0}\nu_j$ corresponding to toric cyclic singularity of order $|c|$. By the same reason, the lattice 
vectors $\nu_{j+1}, \ldots, \nu_n$ define the minimal regular refinement of the cone $\R_{\geq 0} \nu_{j+1} + \R_{\geq 0}\nu_1$
corresponding to toric cyclic singularity of order $|b|$. Therefore, the whole plane $N_\R$ is covered by three cones 
\[ \R_{\geq 0} \nu_1 + \R_{\geq 0}\nu_j, \;\;\R_{\geq 0} \nu_j + \R_{\geq 0}\nu_{j+1}, \;\; \R_{\geq 0} \nu_{j+1} + \R_{\geq 0}\nu_1, \]
both integers $b, c$ must be positive, and lattice verctors $\nu_1, \nu_j, \nu_{j+1}$ span a fan of the singular weighted projective space $\P(1,c,b)$. By above arguments, 
the weighted circular graph $\Gamma_r$ defines the minimal desingularization $\widehat{\P}(1, c,b)$ of the singular weighted projective plane $\P(1,c,b)$. 
Since $\nu_j = (-1, 0)$, $\nu_{j+1} = (0,-1)$, we obtain $\nu_1 =(b,c)$. Note that  the lattice vector $\nu_1$ is contained in the interior of the cone $\sigma(a)$ spanned by $(a, 1)$ and $(a+1,1)$. So we obtain $a < b/c < a+1$, \ie , $a = [\frac{b}{c}]$ and $b>c$. In order to show that $r = a+ \delta = b/c$ we use the isomorphism between the cone 
$\sigma(a)$ and the cone spanned by $(1,0)$ and $(1,1)$ whose primitive interior lattice points we enumerated by 
rational numbers $\delta$. Formally, one could apply induction on the length of path in the tree $\mathcal{OT}(a)$ connecting the root node $\Gamma(a)$ with $\Gamma_r$ by a sequence of combinatorial blow-up. The basis of this induction was explained in \ref{1/2}. It remains to use the correspondence between blow-ups and Farey sums. \hfill $\Box$

\section{Applications to minimal horospherical embeddings}

Now our purpose is to explain an application of toric surfaces $S_r$ 
to the classification of  minimal smooth projective quasi-homogeneous $3$-folds
having a regular  action of the group
${SL}(2)\times{\G}_m$. The complete classificatiion was  obtained by  Z. Guan 
\cite{Gu97,Gu99,Gu00}. His main result in  \cite[Theorem 4.1]{Gu00} 
claims that any smooth minimal projective $3$-fold
having a  quasi-homogeneous regular action of the group
${SL}(2)\times{\G}_m$ belongs to one of the following 5 types:
\begin{enumerate}

\item projective space ${\P}^3$;

\item $3$-dimensional quadric
in $\P^4$;

\item ${\P}^2$-bundles over ${\P}^1$;

\item $\P^1$-bundles over ${\P}^2$,
or over a Hirzebruch surface ${\F}_n$ $(n \geq 0)$;

\item "house models'' $A_{p,q}$ for coprime
positive integers $p,q$ $(p>1)$.

\end{enumerate}

The earlier classification results of Mukai, Umemura and Nakano
about quasi-homogeneous $SL(2)$-manifolds \cite{MU83,Na89} show
that in order to obtain  the Guan's classification
it remains  to consider 
quasi-homogeneous $SL(2)\times \G_m$-threefolds   $X$  such that
a generic $SL(2)$-orbit in $X$ has dimension $2$ (see \cite{Gu97}).
In the investigation of such varieties  
Guan used arguments based on   the Bialynicki-Birula decomposition
for $\G_m$-actions \cite{Bia73} and the Mori theory in dimension $3$
\cite{Mo82}.

This paper suggest an alternative approach to  the 
classification results of Guan  which is based on the fact  that
if a generic $SL(2)$-orbit in $X$ has dimension $2$ 
then such an orbit is a spherical $SL(2)$-variety of rank $1$ and
the $3$-fold $X$ can be considered as a spherical
$SL(2)\times \G_m$-variety of rank $2$.
Therefore, according to the  Luna-Vust theory  \cite{LV83,Kn91},
one can describe the minimal $3$-folds $X$ combinatorially
by means of a $2$-dimensional colored fans. This idea was used in \cite{Ha10}.

We find the language of Luna-Vust especially convenient
for describing the infinite series of minimal ``house models'' $A_{p,q}$ found by Guan.
Note that all minimal $3$-folds  $A_{p,q}$ are  horospherical
$SL(2)\times \G_m$-varieties.  On the other hand, it is known that any horospherical
$SL(2)\times \G_m$-variety is a toric variety (see \cite{Pas08}). In particular, 
all house models $A_{p,q}$ are simultaneously minimal
smooth projective   toric 
$3$-folds which  were described by Guan using $3$-dimensional simple convex polytopes
$\Delta_{p,q}$.  The convex polytopes $\Delta_{p,q}$
has a shape similar to a "house" (this fact probably served as the reason for the use of  the name {\em house model}). 
In general, the classification of minimal smooth 
projective toric $3$-folds is a rather complicated combinatorial
problem. It has been solved only for the
Picard rank $\leq 5$ \cite[Theorem 1.34]{Oda85}.
On the other hand, the minimal toric
$3$-folds $A_{p,q}$ obtained by Guan provide infinite series of 
varieties which may have arbitrary large
Picard rank.

Unfortunately, the definition of  minimal toric $3$-folds
$A_{p,q}$ was given by Guan in a rather  complicated form which
uses sequences of toric blow-ups and blow-downs
relating  a toric $3$-fold $A_{p,q}$ to another toric $3$-fold $A_{p',q'}$
for smaller integers $p',q'$ (see \cite[p.578-580]{Gu99}). 
Smooth toric surfaces $S_r$   with  $r = \frac{p+q}{p}$ considered 
in this paper  provide a 
straightforward   description of  Guan's $3$-folds
$A_{p,q}$.  We illustrate this method in a simplest situation for $p=2$ and $q=1$. 

\begin{exam}
Consider the toric surface $S_{\frac{3}{2}}$ determined by the $2$-dimensional fan 
$\Sigma_{\frac{3}{2}}$.  The contraction of the exceptional curve on $S_{\frac{3}{2}}$ produces a toric surface $S$ given by the fan $\Sigma$.  On the other hand, we create the  colored fan $\Sigma^c_{\frac{3}{2}}$ by attaching color to the vertex $\nu$ with the weight $w=-1$. 
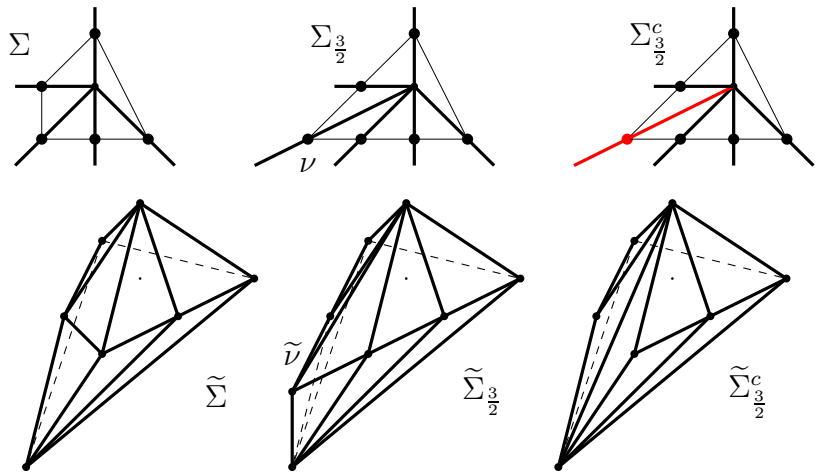
\begin{center}
\begin{tikzpicture}[scale=0.7]
\fill (0,0) circle (2pt);
\fill (0,1) circle (3pt);
\fill (0,-1) circle (3pt);
\fill (-1,-1) circle (3pt);
\fill (-1,0) circle (3pt);
\fill (-2,-1) circle (3pt);
\fill (1,-1) circle (3pt);

\node  at (-2,-1.5) {{$\nu$}}; 

\draw[solid, very thin] (0,1) -- (1,-1) -- (0,-1) -- (-1,-1) -- (-2,-1) -- (0,1);

\draw[solid, very thick] (0,0) -- (0,1.5) ;
\draw[solid, very thick] (0,0) -- (1.5,-1.5) ;
\draw[solid, very thick] (0,0) -- (0,-1.5) ;
\draw[solid, very thick] (0,0) -- (-1.5,-1.5) ;
\draw[solid, very thick] (0,0) -- (-3,-1.5) ;
\draw[solid, very thick] (0,0) -- (-1.5,0) ;

\node [left] at (-1,0.8) {{$\Sigma_{\frac{3}{2}}$}}; 

\begin{scope}[xshift= 6cm]

\fill (0,0) circle (2pt);
\fill (0,1) circle (3pt);
\fill (0,-1) circle (3pt);
\fill (-1,-1) circle (3pt);
\fill (-1,0) circle (3pt);
\fill [color=red] (-2,-1) circle (3pt);
\fill (1,-1) circle (3pt);

\draw[solid, very thin] (0,1) -- (1,-1) -- (0,-1) -- (-1,-1) -- (-2,-1) -- (0,1);

\draw[solid, very thick] (0,0) -- (0,1.5) ;
\draw[solid, very thick] (0,0) -- (1.5,-1.5) ;
\draw[solid, very thick] (0,0) -- (0,-1.5) ;
\draw[solid, very thick] (0,0) -- (-1.5,-1.5) ;
\draw[color=red, solid, very thick] (0,0) -- (-3,-1.5) ;
\draw[solid, very thick] (0,0) -- (-1.5,0) ;

\node [left] at (-1,0.8) {{$\Sigma_{\frac{3}{2}}^c$}}; 
\end{scope}
\begin{scope}[xshift= -6cm]
\node [left] at (-1,0.8) {{$\Sigma$}}; 

\fill (0,0) circle (2pt);
\fill (0,1) circle (3pt);
\fill (0,-1) circle (3pt);
\fill (-1,-1) circle (3pt);
\fill (-1,0) circle (3pt);
\fill (1,-1) circle (3pt);

\draw[solid, very thin] (0,1) -- (1,-1) -- (0,-1) -- (-1,-1) -- (-1,0) -- (0,1);

\draw[solid, very thick] (0,0) -- (0,1.5) ;
\draw[solid, very thick] (0,0) -- (1.5,-1.5) ;
\draw[solid, very thick] (0,0) -- (0,-1.5) ;
\draw[solid, very thick] (0,0) -- (-1.5,-1.5) ;
\draw[solid, very thick] (0,0) -- (-1.5,0) ;
\end{scope}

\end{tikzpicture}

\end{center}

\begin{center}
\begin{tikzpicture}[scale=0.5]
\node  at (2,-3) {{$\widetilde{\Sigma}_{\frac{3}{2}}$}}; 

\fill (0,0) circle (1pt);
\fill (3,0) circle (3pt);
\fill (0,2) circle (3pt);
\fill (-1,1) circle (3pt);
\fill (-1,-2) circle (3pt);
\fill (-2,-1) circle (3pt);
\fill (-3,-3) circle (3pt);
\fill (1,-1) circle (3pt);
\fill (-3,-5) circle (3pt);

\node  at (-3,-2) {{$\widetilde{\nu}$}}; 

\draw[solid, very thick] (-3,-5) -- (3,0) --(0,2) -- (-1, 1) -- (-3,-3) -- (-3,-5); 

\draw[solid, very thick] (-3,-3) -- (3,0) ;
\draw[solid, very thick] (-3,-3) -- (0,2) ;

\draw[solid, very thick] (-1,-2) -- (0,2) ;
\draw[solid, very thick] (1,-1) -- (0,2) ;
\draw[solid, very thick] (-2,-1) -- (0,2) ;

\draw[solid, very thick] (-1,-2) -- (-3,-5) ;
\draw[solid, very thick] (1,-1) -- (-3,-5) ;
\draw[dashed]  (-2,-1) -- (-3,-5) ;
\draw[dashed]  (-1,1) -- (-3,-5) ;

\draw[dashed] (-1,1) -- (3,0);

\begin{scope}[xshift= -7cm]

\node  at (2,-3) {{$\widetilde{\Sigma}$}}; 

\fill (0,0) circle (1pt);
\fill (3,0) circle (3pt);
\fill (0,2) circle (3pt);
\fill (-1,1) circle (3pt);
\fill (-1,-2) circle (3pt);
\fill (-2,-1) circle (3pt);
\fill (1,-1) circle (3pt);
\fill (-3,-5) circle (3pt);

\draw[solid, very thick] (-2,-1) -- (-1,-2) --(3,0) -- (0, 2) -- (-1,1) -- (-2,-1); 

\draw[solid, very thick] (-3,-5) -- (3,0) ;

\draw[solid, very thick] (-1,-2) -- (0,2) ;
\draw[solid, very thick] (1,-1) -- (0,2) ;
\draw[solid, very thick] (-2,-1) -- (0,2) ;

\draw[solid, very thick] (-1,-2) -- (-3,-5) ;
\draw[solid, very thick] (1,-1) -- (-3,-5) ;
\draw[solid, very thick] (-2,-1) -- (-3,-5) ;
\draw[dashed]  (-1,1) -- (-3,-5) ;

\draw[dashed] (-1,1) -- (3,0);

\end{scope}
\begin{scope}[xshift= 7cm]

\fill (0,0) circle (1pt);
\fill (3,0) circle (3pt);
\fill (0,2) circle (3pt);
\fill (-1,1) circle (3pt);
\fill (-1,-2) circle (3pt);
\fill (-2,-1) circle (3pt);
\fill (1,-1) circle (3pt);
\fill (-3,-5) circle (3pt);

\draw[solid, very thick]  (-1,-2) --(3,0) -- (0, 2) -- (-1,1) -- (-2,-1); 

\draw[solid, very thick] (-3,-5) -- (3,0) ;

\draw[solid, very thick] (-1,-2) -- (0,2) ;
\draw[solid, very thick] (1,-1) -- (0,2) ;
\draw[solid, very thick] (-2,-1) -- (0,2) ;

\draw[solid, very thick] (0,2) -- (-3,-5) ;

\draw[solid, very thick] (-1,-2) -- (-3,-5) ;
\draw[solid, very thick] (1,-1) -- (-3,-5) ;
\draw[solid, very thick] (-2,-1) -- (-3,-5) ;
\draw[dashed]  (-1,1) -- (-3,-5) ;

\draw[dashed] (-1,1) -- (3,0);

\node  at (2,-3) {{$\widetilde{\Sigma}^c_{\frac{3}{2}}$}}; 
\end{scope}
\end{tikzpicture}

\end{center}
By considering three $2$-dimensional fans  $\Sigma$,  $\Sigma_{\frac{3}{2}}$ and $\Sigma^c_{\frac{3}{2}}$ as fans of horospherical $SL(2) \times \G_m$-varieties, we obtain 
three  smooth toric $3$-folds $\widetilde{V}$, $\widetilde{V}_{\frac{3}{2}}$ and   $\widetilde{V}^c_{\frac{3}{2}} = A_{2,1}$ which can be also  described by the $3$-dimensional fans $\widetilde{\Sigma}$, $\widetilde{\Sigma}_{\frac{3}{2}}$ 
and $\widetilde{\Sigma}^c_{\frac{3}{2}}$ respectively. The toric fan $\widetilde{\Sigma}_{\frac{3}{2}}$ of the horospherical 
$SL(2) \times \G_m$-variety $\widetilde{V}_{\frac{3}{2}}$ given by the $2$-dimensional uncolored fan $\Sigma_{\frac{3}{2}}$ contains a divisor $D$ isomorphic to $\P^1 \times \P^1$  corresponding to the lattice vector $\widetilde{\nu}$. The divisor $D$  
admits contraction to a curve $\P^1$ in another smooth non-minimal 
toric $3$-fold $\widetilde{V}$ given by the fan $\widetilde{\Sigma}$. On the other hand, there exists another contraction of $D$ 
to $\P^1$ in the minimal toric $3$-fold $\widetilde{V}^c_{\frac{3}{2}}$ given by the fan  $\widetilde{\Sigma}^c_{\frac{3}{2}}$. 

\end{exam}
 
The construction of the minimal house $3$-fold $A_{p,q}$ for arbitrary coprime $p,q$ $p>1$ works analogously. 
We  put $r := \frac{p+q}{p}$, take the complete regular 
fan $\Sigma_r$ of the toric surface $S_r$ and  make  a colored fan 
$\Sigma_r^c$ by attaching the color to the vertex $\nu$ corresponding to the 
exceptional curve. Then $A_{p,q}$ is  the $3$-dimensional 
horospherical $SL(2) \times \G_m$-variety $\widetilde{\Sigma}^c_r$ 
defined by the colored fan $\Sigma_r^c$. 

Now it is easy to describe  the $3$-dimensional
fan of the minimal toric $3$-fold $A_{p,q}$. 
We begin with the $2$-dimensional uncolored fan $\Sigma_r$ $(r =  \frac{p+q}{p})$. 
As a fan of a horospherical $SL(2) \times \G_m$-variety it defines 
a locally trivial
fibration over $\P^1$ whose fiber is isomorphic to a smooth
projective toric surface $S_r$. Let $\nu_1, \ldots, \nu_{n}$ be the primitive lattice
generators of $1$-dimensional cones in the fan $\Sigma_r$ such that $\nu_1$ 
spans the ray $\R_{\geq 0}$ corresponding to the unique exceptional 
curve $C_1 \subset S_r$. Then we have  
$  \nu_2 + \nu_n = \nu_1$. 
Then we put the  $2$-dimensional fan $\Sigma_r$ into $N_\R' \cong \R^3$ with the $\Z$-basis  
$\nu_2, \nu_n, \nu' \in N$. Define another lattice vector $\nu'' \in N'$ by linear relation
\[ \nu'+ \nu'' = \nu_2 + \nu_n  = \nu_1. \]
The $1$-dimensional cones in the fan $\widetilde{\Sigma}_r$ are spanned by $\nu_1, \ldots, \nu_n,\nu',  \nu'' \in N'$
and $3$-dimensional cones in the fan $\widetilde{\Sigma}_r$ are 
\[ \R_{\geq 0}\nu' + \sigma, \;\;  \R_{\geq 0} \nu'' + \sigma \;\; (\sigma \in \Sigma_r). \]
This fan defines the smooth projective toric $3$-fold $\widetilde{V}_r$ which has  locally trivial toric 
fibration over $\P^1$ with the  fiber  $S_r$ because  $ \nu'  +\nu'' = \nu_1 
\in N$. The toric divisor $D$ corresponding to $\nu_1 \in \widetilde{\Sigma}_r$ is isomorphic to the product 
$\P^1 \times \P^1$ and it has two different birational contractions $\psi$, $\varphi$. 

The first contraction 
$\psi\, :\, \widetilde{V}_r \to \widetilde{V}$ does not violate the fibration structure because it contracts $\P^1$ in each fiber to 
a point so that eventually we obtain that $\psi(D) \cong \P^1$ is a section in the fibration over $\P^1$. 
The second contraction $\varphi\, :\, \widetilde{V}_r \to \widetilde{V}^c_r$ does violate the fibration structure because the image 
$\varphi(D) \cong \P^1$ becomes the  indeterminacy locus for the fibration morphism to $\P^1$. This is the reason why 
the variety $\widetilde{V}^c_r$ becomes minimal.  Using results of Brion, we summarize properties of  birational morphisms between smooth horospherical 
$SL(2) \times \G_m$-varieties  in \ref{bir-horo} and \ref{minimal-horo}.

\newpage

\section{Appendix: Spherical and horospherical varieties}
Let $G$ be a connected  reductive algebraic group. A normal
  $G$-variety $X$ is called {\bf quasi-homogeneous} if $X$
contains  a dense open $G$-orbit.  Then $X$ 
is  considered $X$ as a $G$-equivariant
embedding of a dense open $G$-orbit $Gx \cong G/H$, where
$H$ is stabilizer of a point $x \in X$ in the open  $G$-orbit.
A quasi-homogeneous $G$-variety $X$
is called {\bf spherical} if
a Borel subgroup $B \subset G$ contains a dense open
orbit in $X$. If $X$ is a spherical $G$-variety, then the open orbit  $Gx \cong
G/H$ is called a {\bf spherical homogeneous space}, $H$ is
called a {\bf spherical subgroup},  and $X$ is called
a {\bf spherical embedding} of $G/H$.

We will be interested in special spherical subgroups
$H \subset G$ which
contain a maximal unipotent subgroup $U \subset G$.
These subgroups $H$ are called {\bf horospherical} and
equivariant embeddings of $G/H \hookrightarrow X$ are called
{\bf horospherical embeddings}.
The Luna-Vust theory describes  arbitary spherical embeddings  using
some combinatorial objects which are called
{\bf colored fans} (see \cite{LV83}, \cite{Br97},
\cite{Kn91}). Horospherical $G$-equivariant
embeddings  of $G/U$ were considered independently
by Pauer \cite{Pau81,Pau83}. Some further   results on
horospherical embeddings  were contained  in the PhD thesis
of Pasquier \cite{Pas08}. (see also the book of Timashev
\cite{Ti11}).

Let us  recall the Luna-Vust theory for horospherical embeddings $G/H  \hookrightarrow X$.
If   $H \subset G$ is  horospherical subgroup, then the normalizer
$N_G(H)$ is a parabolic subgroup $P \subseteq G$ and
$N_G(H)/H$ is an algebraic torus $T$. Let $M:= {\mathcal X}(T)$
be the lattice of characters
of $T$. We consider $M$  as a sublattice  in the lattices 
of characters  ${\mathcal X}(P) \subseteq {\mathcal X}(B)$ of $P$. Let
$\Pi(G) = \{ \alpha_1, \ldots, \alpha_n\}$ be the set
of simple roots of $G$. We denote by $P_{\alpha}$ the minimal parabolic
subgroup in $G$ corresponding to a simple root $\alpha \in \Pi(G)$.
Then the parabolic subgroup $P = N_G(H)$ is uniquely determined
up to a conjugation by a subset $I := \{ \alpha \in \Pi(G)\; : \; P_\alpha \subseteq
P \}$.
Using the locally trivial fibration  $G/H \to G/P$ (with fibers isomorphic
to $T$),  one obtains a natural bijection between the set of all
irreducible $B$-invariant divisors
${\mathcal D}(G/H):= \{ D_1, \ldots, D_k\}$ in $G/H$ and
the complementary subset of simple roots $\overline{I} := \Pi(G) \setminus I$.
Let $N:= {\rm Hom} (M, \Z)$ be the dual lattice.
There exist a natural map
\[  \rho\; : \; {\mathcal D}(G/H) \to N, \;\; \;\; D_\alpha \mapsto \alpha^\vee|_M \],
where $D_\alpha$  is the  irreducibal $B$-invariant divisor
corresponding to a simple root $\alpha \in \overline{I}$ and $\alpha^\vee|_M$  is
the restriction of the coroot $\alpha^{\vee}\, : \,
{\mathcal X}(B) \to \Z$ to the sublattice $M \subseteq  {\mathcal X}(P)
\subseteq {\mathcal X}(B)$.

Let ${\mathcal G}$ be a finite subset in $N$.
The {\bf rational polyhedral cone} $\sigma = \sigma ({\mathcal G})$
generated by ${\mathcal G}$
is the subset $\sigma \subset N_\R :=
N \otimes_Z \R$  defined as
\[ \sigma:= \{ \sum_{ n \in {\mathcal G} } \R_{\geq 0} n \}. \]
The cone $\sigma$ is called {\bf strictly convex} if
$-\sigma \cap \sigma = 0$.  A {\bf colored cone} is a pair
$(\sigma, {\mathcal F}_\sigma)$
consisting of a strictly convex rational polyhedral cone $\sigma$
and a subset $ {\mathcal F}_\sigma \subseteq D(G/H)$ such that
$\rho(D)$ is a nonzero lattice point in $\sigma$ for all
$D \in  {\mathcal F}_\sigma$.
Every colored cone $(\sigma, {\mathcal F}_\sigma)$ defines
a simple quasi-projective horospherical $G$-embedding
$G/H \hookrightarrow X_{\sigma, {\mathcal F}_\sigma}$ containing 
a unique closed $G$-orbit.
A {\bf face}  of a colored cone  $(\sigma, {\mathcal F}_\sigma)$
is a pair $(\sigma', {\mathcal F}_\sigma')$ such that $\sigma'$
is a face of $\sigma$ and ${\mathcal F}_\sigma'=
\{ D \in {\mathcal F}_\sigma\; : \; \rho(D) \in \sigma' \}$.
A {\bf colored fan} $\Sigma$ is a finite
collection of colored cones such that a face of a colored
cone in $\Sigma$ belongs again to $\Sigma$, and for any
two colored cones  $(\sigma, {\mathcal F}_\sigma)$ and
 $(\sigma', {\mathcal F}_\sigma')$ the pair
 $(\sigma \cap \sigma', {\mathcal F}_\sigma \cap  {\mathcal F}_\sigma')$
is a common face of both colored cones. The Luna-Vust theory
estabishes a bijection between the set of all $G$-equivariant
horospherical embeddings $G/H \hookrightarrow X_\Sigma$ and the set of colored
fans \cite{LV83}.  Many properties of $X_\Sigma$ can be characterized
via combinatorial properties of the colored fan $\Sigma$.
For example, $X_\Sigma$ is {\bf complete} if and only if
\[ N_\R = \bigcup_{(\sigma, {\mathcal F}_\sigma) \in \Sigma} \sigma. \]
There is a natural bijection between colored cones
in $\Sigma$ and $G$-orbits in $X_\Sigma$.
The normal variety $X_\Sigma$ is {\bf locally factorial} if and only if
for every clored cone $(\sigma, {\mathcal F}_\sigma)$ the cone
$\sigma$ is generated by a part of $\Z$-basis of the lattice
$N$ and for any $D \in {\mathcal F}_\sigma$ the lattice vector
$\rho(D)$ belongs to this generating set. The {\bf smoothness criterion}
for horospherical varieties is more subtile than for toric varieties
\cite{Pau83}.

We consider below examples  illustrating   the Luna-Vust theory
in the case $G= SL(2) \times \G_m$ and $H= U \times 1$,
where $U$ is a maximal unipotent subgroup in $G$ consisting of
upper triangular matrices with eigenvalues $1$.
The normalizer $N_G(H)$ is a Borel subgroup $B$ in $G$:
\[ B := \left\{ \left( \begin{pmatrix} \lambda & * \\ 0 &
\lambda^{-1} \end{pmatrix}, \mu \right)\;\; :  \;\; \lambda, \mu \in \C^*
\right\}.
\]
Then $\{\lambda, \mu \}$ can be used to define a natural basis of the lattice
of characters $M$ of $B$ so that we can identify this group
with $\Z^2$. There exists exactly one irreducible
$B$-invariant divisor $D$ in $G/H$. Using the dual basis
in $N = {\rm Hom}(M, \Z)$ we can identify
 $\rho(D) \in N$ with the lattice point $(1,0)$.
The quotient $G/H$ is a quasi-affine variety
$(\C^2 \setminus (0,0))\times \C^*$.

\begin{exam}
Consider the natural simple $G$-embedding $G/H \hookrightarrow  \C^3$
consisting of $4$ $G$-orbits: $(0,0,0)$, $(0,0) \times \C^*$,
$(\C^2 \setminus \{ (0,0) \}) \times 0$ and the dense $G$-orbit.
The embedding $G/H \hookrightarrow  \C^3$ corresponds to the colored
fan $\Sigma_1$ consisting of the $2$-dimensional colored cone $(\R_{\geq 0}(1,0) +
\R_{\geq 0}(0,1), D)$ and
all its faces. 
\end{exam}

\begin{exam}
Another natural 
$G$-embedding of $G/H \hookrightarrow  \P^3$ extends the
previous one by adding at infinity the projective plane $\P^2$ consisiting
of $3$ more $G$-orbits: a point, $\C^2 \setminus (0,0)$, and
$\P^1$. This embedding of $G/H$ corresponds to a colored
fan $\Sigma_2$ having  additional $3$ colored cones
$(\R_{\geq 0}(1,0) + \R_{\geq 0}(-1,-1), D)$,
$(\R_{\geq 0}(-1,-1), \emptyset)$, and  $(\R_{\geq 0}(0,1) +
\R_{\geq 0}(-1,-1), \emptyset).$
\end{exam}

There is a natural operation which produces from any colored fan $\Sigma$
the  uncolored fan  $\Sigma'$ by replacing every colored cone
$(\sigma, {\mathcal F}_\sigma)$ by the pair $(\sigma, \emptyset)$.
The bijection
between cones in $\Sigma'$ and $\Sigma$  defines
a proper surjective $G$-equivariant morphism $X_{\Sigma'} \to X_\Sigma$.
Let $Y_{\Sigma'}$ be $T$-equivariant embedding of the torus $N_G(H)/H$
corresponding to  the uncolored
fan $\Sigma'$. Then the horospherical embedding $G/H \hookrightarrow
X_{\Sigma'}$ can be described as a toric fibration over $G/P$ with the
fiber $Y_{\Sigma'}$.

\begin{exam} 
If we apply the above procedure   to the colored fan $\Sigma_2$
corresponding to the horospherical embedding $G/H
\hookrightarrow \P^3$ we obtain the uncolored fan $\Sigma'_2$ defining
the horospherical embedding
$$G/H \hookrightarrow
X_{\Sigma'_2} \cong
\P_{\P^1}({\mathcal O} \oplus {\mathcal O}  \oplus {\mathcal O}(1)).$$
The horospheric variety $X_{\Sigma'_2}$ is a $\P^2$-bundle over $G/B \cong \P^1$
and it
can be obtained as a blow-up
of the $G$-invariant line  in $\P^3 = X_{\Sigma_2}$ which is the closure
of the $1$-dimensional $G$-orbit corresponding to the colored
cone $(\R_{\geq 0} (1,0), D)$.
\end{exam}

Equivariant birational contractions in the Mori theory
for spherical varieties were  described  by Brion \cite{Br93}
using methods similar to the ones for toric varieties
\cite{Re83}. It is important to remark  that  for  spherical varieties
extremal rays in the cone of effective $1$-cycles are generated
by classes of closures of $1$-dimensional $G$-orbits. Therefore,
the corresponding  Mori contractions are  $G$-equivariant.

In contrast to toric varieties the $G$-equivariant birational geometry
of spherical varieties is more rich due to existence of colors.
For smooth projective spherical varieties of rank $2$ Brion has
generalized the factorization theorem for equivariant
birational maps between smooth projective  toric surfaces \cite{Br94}.

For our purposes, we summarize some known  results in the case
of $3$-dimensional horospherical $SL(2) \times \G_m$-varieties:

\begin{theorem} Let $X_\Sigma$ be a
projective $3$-dimensional  horospherical
$SL(2) \times \G_m$-variety corresponding to a  $2$-dimensional colored fan
$\Sigma$. Then $X_\Sigma$ is smooth if and only if two conditions
are safisfied:

\begin{enumerate}

\item every cone $\sigma$ in $\Sigma$ is generated by a part of
a $\Z$-basis of the lattice $N$;

\item  if a cone $(\sigma, D) \in \Sigma$
has the color $D$, then the lattice point $\rho(D) \in N$ must be
one of its primitive  generators.
\end{enumerate}
\end{theorem}

\begin{theorem} \label{bir-horo} Let $X_\Sigma$ be a smooth
projective $3$-dimensional  horospherical
$SL(2) \times \G_m$-variety corresponding to a $2$-dimensional
colored fan
$\Sigma$. Then there exist exactly $3$ types of equivariant blow-ups of
$X_\Sigma$:
\begin{enumerate}

\item the blow up of a $SL(2) \times \G_m$-fixed point in
$X_\Sigma$ corresponding to a $2$-dimensional colored cone
$(\sigma, D)$, where $\sigma = \R_{\geq 0} \nu + \R_{\geq 0} \nu'$ for
some $\Z$-basis $\nu, \nu'$ of $N$. If $\nu' = \rho (D)$, then we set
$\nu'':= \nu + \nu'$ and  the colored
fan of the blow-up
is obtained by the  subdivision of $(\sigma, D)$ into two $2$-dimensional
subcones $(\sigma', \emptyset)$ and   $(\sigma'', D)$, where
 $\sigma' = \R_{\geq 0} \nu + \R_{\geq 0} \nu''$ and  $\sigma'' = \R_{\geq 0} \nu' +
\R_{\geq 0} \nu''$.

\item the blow up of a closed $1$-dimensional $SL(2) \times \G_m$-orbit
in $X_\Sigma$ corresponding to a $2$-dimensional uncolored cone
$(\sigma, \emptyset)$, where $\sigma = \R_{\geq 0} \nu + \R_{\geq 0} \nu'$ for
some $\Z$-basis $\nu, \nu'$ of $N$. Again  we set
$\nu'':= \nu + \nu'$ and  the colored
fan of this  blow-up
is obtained by the  subdivision of $(\sigma, \emptyset)$
into two $2$-dimensional
uncolored subcones $(\sigma', \emptyset)$ and   $(\sigma'', \emptyset)$,
where  $\sigma' = \R_{\geq 0} \nu + \R_{\geq 0} \nu''$ and
$\sigma'' = \R_{\geq 0} \nu' + \R_{\geq 0} \nu''$.

\item the blow up of the closure of a  $1$-dimensional
$SL(2) \times \G_m$-orbit
in $X_\Sigma$ corresponding to the $1$-dimensional colored cone
$(\R_{\geq 0} \rho(D), D)$. The fan of this blow-up
is obtained by deleting the color $D$ from all colored
cones of $\Sigma$.
\end{enumerate}
\end{theorem}

\begin{coro}
\label{minimal-horo}
A smooth
projective $3$-dimensional  horospherical
$SL(2) \times \G_m$-variety defined by a colored fan $\Sigma$
is minimal (i.e. does not admit any Mori contraction to another
smooth projective $3$-fold)  if and only if $\Sigma$ does not contain
two $2$-dimensional colored cones $( \R_{\geq 0} \nu + \R_{\geq 0} \nu'',
{\mathcal F}_1)$ and  $( \R_{\geq 0} \nu'' + \R_{\geq 0} \nu',
{\mathcal F}_2)$ with the uncolored common face  $(\R_{\geq 0} \nu'', \emptyset)$
such that $\nu'' = \nu + \nu'$.
\end{coro}

\end{document}